\documentclass[12pt,diff_eng]{article}
\usepackage[english]{babel}
\usepackage{amssymb,amscd,amsmath}
\newtheorem{theorem}{Theorem}

\newtheorem{corollary}{Corollary}
\newcommand{\be}{\begin{equation}}
\newcommand{\ee}{\end{equation}}
\newcommand{\bea}{\begin{eqnarray}}
\newcommand{\eea}{\end{eqnarray}}
\newcommand{\nn}{\nonumber}
\newcommand{\rd}{\partial}
\newcommand{\mc}{{\mathbb C}}

\begin{document}


 \begin{center}
    \font\titlerm=cmr10 scaled\magstep4
    \font\titlei=cmmi10 scaled\magstep4
    \font\titleis=cmmi7 scaled\magstep4
     \centerline{\titlerm Green's Function For Linear Differential}
      \vspace{0.5cm}
     \centerline{\titlerm  Operators In One Variable }

    \vspace{1.5cm}
    \noindent{{
        Adel Kassaian{$^*$}\footnote{E-mail:
        a.kassaian@gmail.com},
         }}\\

  \end{center}

\begin{abstract}

General formula for causal Green's function of linear differential
operator of given degree in one variable, $
({\rd^{\,n}_{x}}+\sum^{n-1}_{k=0} P_{k}(x) {\rd^{k}_{x}}) $, is
given according to coefficient functions of differential operator as
a series of integrals. The solution also provides analytic formula
for fundamental solutions of corresponding homogenous linear
differential equation, $\big({\rd^{\,n}_{x}}+\sum^{n-1}_{k=0}
P_{k}(x) {\rd^{k}_{x}}\big)\,y(x)=0$, as series of integrals.
Furthermore, multiplicative property of causal Green's functions is
shown and by which explicit formulas for causal Green's functions of
some classes of decomposable linear differential operators are
given. A method to find Green's function of general linear
differential operator of given degree in one variable with arbitrary
boundary condition according to coefficient functions of
differential operator is demonstrated.
\end{abstract}

\section{Converting initial value problem for ordinary linear differential equation into  Volterra's integral equation}
Initial value problem for linear differential equation of degree $n$
in one variable, \be\big({\rd^{\,n}_{x}}+\sum^{n-1}_{k=0} P_{k}(x)
{\rd^{k}_{x}}\big)\,y(x)=g(x),\label{yak}\ee  can be converted to
Volterra's integral equation of second kind. Following the reference
\cite{ik}, the key relation in the procedure is the identity;

 \be
 \int_{a}^{x}dz_1
\int_{a}^{z_{1}} dz_{2} \cdots\int_{a}^{z_{r-1}}\,\,dz_{r}\,\,
F(z_{r})= \int_{a}^{x} dz {(x-z)^{r-1}\over (r-1)!}F(z), \nn\ee
which can be easily proved for arbitrary function $F(x)$ using
integration by part and induction. For initial condition
$\rd^{\,i}_{x}\,y(a)=c_i$ for $ i=0,1,..,n-1$, putting $u(x)=
\rd_{x}^{\,n}\,y(x)$ in the above relation and setting $r=n-k$, for
$k=0,...,n-1$ we have;

\be \rd_{x}^{\,k}\,y(x)-\sum_{i=k}^{n-1}
c_{i}{(x-a)^{i-k}\over(i-k)!}=\int_{a}^{x} dz {(x-z)^{n-k-1}\over
(n-k-1)!}\,u(z).\label{head} \ee

Inserting $\rd_{x}^{\,n}\,y(x)=u(x)$ and
$\rd_{x}^{\,k}\,y(x)=\sum_{i=k}^{n-1}
c_{i}{(x-a)^{i-k}\over(i-k)!}+\int_{a}^{x} dz {(x-z)^{n-k-1}\over
(n-k-1)!}\,u(z)$ for $k=1,...,n-1$ into the left hand side of
equation (\ref{yak}) we get; $
 \big({\rd^{\,n}_{x}}+\sum^{n-1}_{k=0} P_{k}(x)
{\rd^{k}_{x}}\big)\,y(x)=u(x)+\sum_{k=0}^{n-1}\sum_{i=k}^{n-1}
c_{i}P_{k}(x){(x-a)^{i-k}\over(i-k)!}+\int_{a}^{x} dz
\sum^{n-1}_{k=0}P_{k}(x){(x-z)^{i-k}\over (i-k)!}\,u(z)$. Therefore
we have the following Volterra's equation for $u(x)$;

\be u(x)+\int_{a}^{x} dz\,
  K(x,z) \,u(z)=g(x)+S(x),\label{volttera} \ee

where
$K(x,z)=\big(\sum_{k=0}^{n-1}P_{k}(x){{(x-z)^{n-k-1}}\over(n-k-1)!}\big)$
and $S(x)=-\sum_{k=0}^{n-1}\sum_{i=k}^{n-1}
c_{i}P_{k}(x){(x-a)^{i-k}\over(i-k)!}\,=\,
-\sum_{i=0}^{n-1}\sum^{i}_{k=0}
c_{i}{{P_{k}(x)(x-a)^{i-k}}\over(i-k)!}$. By setting $k=0$ in
equation (\ref{head}), the solution $y(x)$ is given by;

\be y(x)=D(x)+\int_{a}^{x} dz {(x-z)^{n-1}\over (n-1)!} u(z).\ee

where $D(x)=\sum_{i=0}^{n-1} c_{i}{(x-a)^{i}\over i !}$.

\section{Casual Green's function for linear differential
operators in one variable}
 \label{tecneq}

For equation (\ref{yak}) with initial condition
$\rd^{\,i}_{x}y(a)=0$ for $ i=0,1,..,n-1$, the corresponding
Volterra's equation is given by;

\be u(x)+\int_{a}^{x} dz\,
  K(x,z) \,u(z)=g(x),\nn \ee

where
$K(x,z)=\big(\sum_{k=0}^{n-1}P_{k}(x){{(x-z)^{n-k-1}}\over(n-k-1)!}\big)$
and $y(x)=\int_{a}^{x} dz {(x-z)^{n-1}\over (n-1)!} u(z)$.For
$g(x)\in L^2[a,b]$, the condition $(\int^{b}_{a}\int^{b}_{a} dx\,
dy\, |K(x,y)|^2)<\infty $ is sufficient condition for existence of
unique solution in $L^2[a,b]$, given by iteration (e.g. see
\cite{it}). Clearly this conditions can be satisfied if $P_{i}(x)$
(\,$i=0,1,...,n-1$) and $g(x)$ functions are taken to be continuous
on $[a,b]$. Therefore we can state the following theorem;

\begin{theorem}\label{FinalGreen's}{The Green's function for inhomogeneous linear differential equation
 $\big({\rd^{\,n}_{x}}+\sum^{n-1}_{k=0} P_{k}(x)\,{\rd^{k}_{x}}
 \big)\,y(x)=g(x),$
  where $P_{i}(x)$ \,(\,$i=0,1,...,n-1$) and $g(x)$ are in $\mc[a,b]$,
  with the boundary condition; $\rd^{\,i}_{x}y(a)=0$ for $i=0,1,..,n-1$, is given
  by;
 \be G(x,y)=\theta(x-y)
\Big({(x-y)^{n-1}\over(n-1)!}+\int^{x}_{y}dz\,{(x-z)^{n-1}\over(n-1)!}\,R(z,y)\Big),\label{maryam}\ee
where \bea R(x,y)=h(x,y)+\sum_{r=2}^{\infty}\int_{y}^{x}dz_1
\int_{y}^{z_{1}} dz_{2} \cdots\int_{y}^{z_{r-2}}\,\,dz_{r-1}\,
h(x,z_1)\,h(z_{1}, z_{2})\nn\\\cdots h(z_{r-1}, y),\eea  and
$\displaystyle\,\,
h(x,y)=-\sum_{k=0}^{n-1}P_{k}(x){{(x-y)^{n-k-1}}\over(n-k-1)!}.$ The
solution to inhomogeneous linear differential equation (\ref{yak}})
for $x\in[a,b]$ is then given by \,\,$ y(x)=\int_{a}^{\infty} dz
\,G(x,z) \,g(z).$ \\
\end{theorem}

{\bf{Proof}}. In order to prove (\ref{maryam}) is Green's function
of (\ref{yak}) it is enough to prove that for the two variables
function; \be
T(x,y)=\Big({(x-y)^{n-1}\over(n-1)!}+\int^{x}_{y}dz\,{(x-z)^{n-1}\over(n-1)!}\,R(z,y)\Big),\label{hasan}\ee
we have $\big({\rd^{\,n}_{x}}+\sum^{n-1}_{k=0} P_{k}(x)
{\rd^{k}_{x}} \big)\,T(x,y)=0$ and $(\rd^{\,i}_{x} T(x,y))|_{x=y}=0$
for $i=0,1,..,n-2$ and $(\rd^{\,n-1}_{x} T(x,y))|_{x=y}=1$ [e.g. see
\cite{ij}]. This can be easily done by noting; \be \rd^{\,i}_{x}
T(x,y)=\Big({(x-y)^{n-i-1}\over(n-i-1)!}+\int^{x}_{y}dz\,{(x-z)^{n-i-1}\over(n-i-1)!}\,R(z,y)\Big)
\hspace{0.5cm} i=0,1,..,n-1\label{magic}\ee \be \rd^{\,n}_{x}
T(x,y)= R(x,y), \nn\ee and therefore, \bea
({\rd^{\,n}_{x}}+\sum^{n-1}_{k=0} P_{k}(x) {\rd^{k}_{x}}
)\,\,T(x,y)=&&R(x,y)+\sum^{n-1}_{k=0}
{ P_{k}(x)(x-y)^{n-k-1}\over(n-k-1)!}\nn\\&&+\sum^{n-1}_{k=0}P_{k}(x)\int^{x}_{y}dz\,{(x-z)^{n-k-1}\over(n-k-1)!}R(z,y)\nn\\
=&& R(x,y)-h(x,y)-\int^{x}_{y} dz\, h(x,z)\,R(z,y) \nn\\=&&
R(x,y)-h(x,y)-\big(R(x,y)-h(x,y)\big)=0.\nn\eea

In the last line we used $\int^{x}_{y} dz\,
h(x,z)\,R(z,y)=\big(R(x,y)-h(x,y)\big)$, which comes from definition
of $R(x,y)$. By using (\ref{magic}) we have $(\rd^{\,i}_{x}
T(x,y))|_{x=y}=0$ for $i=0,1,..,n-2$ and $(\rd^{\,n-1}_{x}
T(x,y))|_{x=y}=1$.\\

The Green's function for (\ref{yak}) with mentioned boundary
condition is called causal solution which by method of variation of
parameters is given by;

\be G(x,y)= \big(\sum_{i=1}^{n}{W_{i}(y)\, u_{i}(x) \over W(y)}
\big) \theta(x-y), \label{korner}\ee

where $u_{1}(x), u_{2}(x),...,u_{n}(x)$ are fundamental solutions of
corresponding homogeneous differential equation;
\,$({\rd^{\,n}_{x}}+\sum^{n-1}_{k=0} P_{k}(x) {\rd^{k}_{x}}
 )\,u_{i}(x)=0$. $W(y)$ is the
Wronskian and $W_{i}(y)$ is the Wronskian with its $i^{\,th}$ column
in determinant is replace by $(0,0,..,0,1)$. Comparing this result
with (\ref{maryam}) we have the identity;

\be \sum_{i=1}^{n}{W_{i}(y)\, u_{i}(x) \over W(y)}
\,=\,{(x-y)^{n-1}\over(n-1)!}+\int^{x}_{y}dz\,{(x-z)^{n-1}\over(n-1)!}\,R(z,y).\label{kasai}\ee

For linear differential operator of first degree ($n$=$1$), like
$\rd_{x}-P(x)$, the causal Green's function using [theorem
\ref{FinalGreen's}] is equal to; \bea(\rd_{x}-P(x))^{-1}&&=
\theta(x-y)\Big(1+\sum_{k=1}^{\infty}\int_{y}^{x}dz_1
\cdots\int_{y}^{z_{n-2}}dz_{k-1}\int_{y}^{z_{k-1}}dz_{k}\,
P(z_1)\cdots P(z_{k-1})P(z_{k})\Big)\nn\\&&=
\theta(x-y)\Big(1+\sum_{k=1}^{\infty}{1\over
k!}\int_{y}^{x}\cdots\int_{y}^{x}\int_{y}^{x}dz_1 \cdots
dz_{k-1}dz_{k}\, P(z_1) \cdots P(z_{k-1})P(z_{k})\Big)\nn
\\&&= \theta(x-y)e^{\int^{x}_{y}dz P(z)}.\nn\eea  For linear
differential operator of degree two in the form of;
$(\rd^{2}_{x}-P(x))$, the causal Green's function by using [theorem
\ref{FinalGreen's}] is given by $T_{s}(x,y)\theta(x-y)$ where; \bea
T_{s}(x,y)=\Big\{(x-y)+\sum_{k=1}^{\infty}\Big(\int_{y}^{x}dz_1
\cdots\int_{y}^{z_{k-2}}dz_{k-1}\int_{y}^{z_{k-1}}dz_{k}\nn\\\,
(x-z_1)P(z_1)(z_1-z_2)P(z_2)(z_2-z_3)\cdots
(z_{k-1}-z_{k})P(z_{k})(z_{k}-y)\Big)\Big\}.\label{shahla}\eea For
example
$(\rd^{2}_{x}-x)^{-1}=\theta(x-y)\,(x-y)+\theta(x-y)\int^{x}_{y}dz(x-z)z(z-y)
+\theta(x-y)\int^{x}_{y} dt\int^{t}_{y} dz\, \,\big((x-t)t
(t-z)z(z-y)\big)+
 \cdots = \theta(x-y) \Big(\big(x  - y\big)+ \big({x^4 \over 12} - {(x^3 y) \over 6}
 + {(x y^3) \over 6} - {y^4 \over 12}\big)+ \big({x^7 \over 504}
  - {(x^6 y) \over 180}+ {(
 x^4 y^3) \over 72} -
  {(x^3 y^4) \over 72} + {(x y^6) \over 180} -
  {y^7 \over 504}\big)+\cdots \Big),
$ which is consistent with solution
$\displaystyle(\rd^{\,2}_{x}-x)^{-1}=\theta(x-y)\big({-\mathrm{Ai}(x)\mathrm{Bi}(y)
 + \mathrm{Ai}(y)\mathrm{Bi}(x)\over{\mathrm{Ai}(y)\mathrm{Bi}'(y)
 - \mathrm{Ai}'(y)\mathrm{Bi}(y)}}\big)$ derived by
 (\ref{korner}).\\

It can be seen from (\ref{maryam}) that if $P_{i}(x)$
(\,$i=0,1,...,n-1$) functions are smooth on
 $[a,b]$ then $T(x,y)$, given by (\ref{hasan}), is smooth function on
 $[a,b]\times[a,b]$, in which case we state the following theorem;

\begin{theorem}\label{adel}{If $T_{1}(x,y)\theta(x-y)$ and $T_{2}(x,y)\theta(x-y)$ are
causal Green's functions for linear differential operators
$\mathcal{O}_{1}(x,\rd_{x})=\big({\rd^{\,n}_{x}}+\sum^{n-1}_{k=0}
P_{k}(x) {\rd^{k}_{x}} \big)$ and
$\mathcal{O}_{2}(x,\rd_{x})=\big({\rd^{\,m}_{x}}+\sum^{m-1}_{k=0}
q_{k}(x) {\rd^{k}_{x}} \big)$ respectively (assuming $P_{i}(x)$'s
and $q_{i}(x)$'s functions are in ${\mc}^{\infty}[a,b]$) then
$T_{3}(x,y)\theta(x-y)$ where, \be T_{3}(x,y)=\int_{y}^{x} dz\,
T_{2}(x,z)\, T_{1}(z,y), \ee is the causal Green's function for
linear differential operator
$\mathcal{O}_{3}(x,\rd_{x})=\mathcal{O}_{1}(x,\rd_{x}).\mathcal{O}_{2}(x,\rd_{x})$}
\end{theorem}

{\bf{Proof.}} By assumption;
$\mathcal{O}_{1}(x,\rd_{x})\,T_{1}(x,y)=0$ and $(\rd^{\,i}_{x}
T_{1}(x,y))|_{x=y}=0$ for $i=0,1,..,n-2$ and $(\rd^{\,n-1}_{x}
T_{1}(x,y))|_{x=y}=1$ also
$\mathcal{O}_{2}(x,\rd_{x})\,T_{2}(x,y)=0$ and $(\rd^{\,i}_{x}
T_{2}(x,y))|_{x=y}=0$ for $i=0,1,..,m-2$ and $(\rd^{\,m-1}_{x}
T_{2}(x,y))|_{x=y}=1$, therefore we have; \be \rd^{\,i}_{x}
\,T_{3}(x,y)=\int_{y}^{x} dz\,(\rd^{\,i}_{x}( T_{2}(x,z)))\,
T_{1}(z,y), \hspace{0.5cm} i=0,1,..,m-1\hspace{0.5cm}\label{kant}\ee
\be \rd^{\,m}_{x} T_{3}(x,y)= T_{1}(x,y)+\int_{y}^{x}
dz\,(\rd^{\,m}_{x}( T_{2}(x,z)))\, T_{1}(z,y), \label{hapar}\ee \bea
\rd^{\,i}_{x} T_{3}(x,y)=
\rd^{\,i-m}_{x}(T_{1}(x,y))+\rd^{\,i-m}_{x}\big(\int_{y}^{x}
dz\,&(\rd^{\,m}_{x}( T_{2}(x,z)))\, T_{1}(z,y)\big),
\nn\\&i=m+1,..,m+n-1.\label{hava}\eea Concentrating on the second
term in (\ref{hava}), we have for $k=1,2,..,n-1$;

\bea \rd^{\,k}_{x}\big(&\int_{y}^{x}
dz\,(\rd^{\,m}_{x}(T_{2}(x,z)))\, T_{1}(z,y) \big)=
\Big\{\sum^{k-1}_{j=0}\rd_{x}^{j}\Big((\rd^{\,m+k-1-j}_{x}
T_{2}(x,z))|_{z=x} \nn\\&  T_{1}(x,y) \Big)\Big\}+\big(\int_{y}^{x}
dz\,(\rd^{\,m+k}_{x}( T_{2}(x,z)))\, T_{1}(z,y)
\big)\nn\\&=\{\sum^{k-1}_{j=0}\sum^{j}_{r=0}{j \choose
r}\Big((\rd^{\,m+k-1-j+r}_{x}
T_{2}(x,z))|_{z=x}\,\rd_{x}^{j-r}T_{1}(x,y) \Big)\}\nn\\
&+\big(\int_{y}^{x} dz\,(\rd^{\,m+k}_{x}( T_{2}(x,z)))\, T_{1}(z,y)
\big)\label{shadi}\eea

From (\ref{kant}),(\ref{hapar}), (\ref{hava}) and (\ref{shadi}) we
have $(\rd^{\,i}_{x} T_{3}(x,y))|_{x=y}=0$ for
$i=0,1,...,m+n-3,m+n-2$ and $(\rd^{\,m+n-1}_{x}
T_{3}(x,y))|_{x=y}=1$. On the other hand; \bea
\mathcal{O}_{2}(x,\rd_{x})
T_{3}(x,y)&=&\big({\rd^{\,m}_{x}}+\sum^{m-1}_{k=0} q_{k}(x)
{\rd^{\,k}_{x}} \big) \int_{y}^{x} dz\, T_{2}(x,z)\,
T_{1}(z,y)\nn\\&=&\rd_{x}\big(\int_{y}^{x} dz\, \rd^{\,m-1}_{x}
T_{2}(x,z)\, T_{1}(z,y)\big)\nn\\&&+ \int_{y}^{x}
dz\,\big(\sum^{m-1}_{k=0} q_{k}(x) {\rd^{\,k}_{x}} T_{2}(x,z)\big)\,
T_{1}(z,y)\nn\\&=& T_{1}(x,y)+\int^{y}_{x}
dz\,\mathcal{O}_{2}(x,\rd_{x}) T_{2}(x,z)
T_{1}(z,y)\nn\\&=&T_{1}(x,y).\nn\eea Therefore
$\mathcal{O}_{1}(x,\rd_{x}).
\mathcal{O}_{2}(x,\rd_{x})T_{3}(x,y)=\mathcal{O}_{1}(x,\rd_{x})T_{1}(x,y)=0.$
\\

 The following corollary comes as a consequence;

\begin{corollary}\label{lema786}   {Causal Green's function for
differential operator,\be {\mathcal
 O}(x,\rd_{x})=\big(\rd_{x}-p_{1}(x)\big)\,\big(\rd_{x}-p_{2}(x)\big)
\,\cdots\big(\rd_{x}-p_{n}(x)\big),
\label{formallimit2}\hspace{1cm}\ee where $p_{\,i}(x)\in
{\mc}^{\infty}[a,b]$ (for $i=1,...,n$)is given by; \bea
G(x,y)=\theta(x-y)&\int_{y}^{x}dz_1 \int_{y}^{z_{1}}dz_{2}
\cdots\int_{y}^{z_{r-2}}dz_{n-1}\big(\,{e^{\int^{x}_{z_{1}}dt_{n}
\,p_{n}(t_{n})}}\nn\\&{e^{\int^{z_{1}}_{z_{2}}dt_{n-1}
\,p_{n-1}(t_{n-1})}}\cdots{e^{\int^{z_{n-1}}_{y} dt_{1}
\,p_{1}(t_{1})}}\big)\label{iraj1}\eea}
\end{corollary}

For example for differential operator
$\mathcal{O}(x,\rd_{x})=\rd^{2}_{x}+3x\rd_{x}+(2x^{2}+2)$, since
$(\rd_{x}+x)(\rd_{x}+2 x)=\rd^{2}_{x}+3x\rd_{x}+(2x^{2}+2)$, by
using result  (\ref{iraj1}) one gets $G(x,y) =\sqrt{\pi\over2}\ \,
e^{\,y^{2}-{x^{2}\over 2}}\, \{\mbox{Erf}
({{x\over\sqrt2}})-(\mbox{Erf}({{y\over\sqrt2}}))\}\,\,\theta(x-y)$.\\

\begin{corollary}\label{hasti}  {\, Causal Green's function for linear differential operator ; \be{\mathcal
 O}(x,\rd_{x})=\sum_{k=0}^{n}\alpha_{k}{\rd^{\,k}_{x}}\, ,\ee where
 $\alpha_{k}\in \mc$ and $\alpha_{n}\ne 0$, is given by;
\be G(x,y)\,\,=\,{\theta(x-y)\over\alpha_{n}}\int_{y}^{x}dz_1
\,\int_{y}^{z_1}\,dz_2 \cdots \int_{y}^{z_{n-2}}dz_{n-1}\,\,
e^{\big( \beta_{1}(x-z_1)+\beta_{2}(z_1-z_2)\cdots
+\beta_{n}(z_{n-1}-y)\big)},\ee where $\beta_{1}$, $\beta_{2}\cdots
\beta_{n}$ are $n$ complex roots of equation
$\displaystyle\sum_{i=0}^{n}\alpha_{i} {X}^{\,i}=0$}.
\end{corollary}

{\bf Proof}. Differential operator ${\mathcal
 O}(x,\rd_{x})=\sum_{i=0}^{n}\alpha_{i}{\rd^{\,i}_{x}}$, according to {\em{Fundamental theorem of
 algebra}},
can be written as, $\,
\sum_{i=0}^{n}\alpha_{i}\rd_x^{{\,i}}=\alpha_n
(\rd_x-\beta_1)(\rd_x-\beta_2)\cdots(\rd_x-\beta_n)$.
 Therefore by using (\ref{iraj1}) the result is proved. \\

For example $ (\rd^{2}_{x}-\omega^2)^{-1} = \theta(x-y)\int_{y}^{x}
dz_{1}
e^{(\omega(x-z_1)-\omega(z_1-y))}={{\sinh\omega(x-y)}\over\omega}\theta(x-y)$\,\,\
 and also $ (\rd^{3}_{x}-i
\alpha\rd^{2}_{x}-\omega^2\rd_{x}+i\alpha\omega^2)^{-1} =
\theta(x-y)\int_{y}^{x} dz_{1}(\int_{y}^{z_1} dz_{2}
e^{(\omega(x-z_1)-\omega(z_1-z_2)+i \alpha(z_2-y))}) $ $
\hspace{0.5cm}=\theta(x-y)\big({{e^{\,\omega(x-y)}-e^{i\alpha(x-y)}}\over{\alpha^2+\omega^2}}
-{\sinh[\omega(x-y)]\over{i\,\alpha\,\omega+\omega^2}}\big). $\\

Lets consider a differential operator in form of \be
\mathcal{O}(x,\rd_x)=-\rd_{x}^{2}+v(x).\label{hamed}\ee By
decomposing it into two firs degree differential operators;
$-(\rd_{x}^{2}-v(x))=-(\rd_{x}-p(x))(\rd_{x}-q(x))$, we have
consequently $q(x)=-p(x)$ and $ p(x)^2-\rd_{x}p(x)=v(x)$. Therefore
according to (\ref{iraj1}) the causal the Green's function is given
by; \be \big(-\rd_{x}^{2}+v(x)\big)^{-1}=-\theta(x-y)\int^{x}_{y} dz
e^{(-\int^{x}_{z}dt p(t)+\int^{z}_{y}dt' p(t'))},\label{sarvin}\ee

where  $p(x)$ is solution for first order nonlinear differential
equation $ p \,(x)^2-\rd_{x} p\, (x)=v(x)$. This is just Riccati
equation, thus the answers to  $ p \,(x)^2-\rd_{x} p\, (x)=v(x)$ are
given by solutions of homogenous differential equation $
(-\rd_{x}^{2}+v(x))u_{1,2}(x)=0 $ where $p(x)=-({u_{1,2}'\over
u_{1,2}})$. Inserting $p(x)=-({u_{1}'\over u_{1}})$ into solution
(\ref{sarvin}), we have; $
\big(-\rd_{x}^{2}+v(x)\big)^{-1}=-\theta(x-y)u_{1}(x)u_{1}(y)\int^{x}_{y}
dz ({1\over {u(z)}^2}).$ Considering the relation $ u_{2}(z)=
u_{1}(z)\int dz {1\over {u_{1}(z)}^{\,2}}$ (valid for homogenous
differential equation $ (-\rd_{x}^{2}+v(x))u_{1,2}(x)=0 $) the
solution (\ref{sarvin}) becomes the standard solution,
$\big(-\rd_{x}^{2}+v(x)\big)^{-1}=-\theta(x-y)\Big(u_{2}(x)u_{1}(y)-u_{1}(x)u_{2}(y)\Big).$
\\

Considering [theorem \ref{adel}] one can introduce the following
infinite non-abelian group of operators on a subspace of
${\mc}^{\infty}[a,b]$. We call it "Lalescu Group";

\begin{itemize}{\item{\bf Lalescu Group.}} {\em The Group of
differential operators of the form;
$\big({\rd^{\,n}_{x}}+\sum^{n-1}_{k=0} P_{k}(x) {\rd^{k}_{x}}\big)$
of all finite order, $n\geq 0$, where $P_{k}(x)\in
{\mc}^{\infty}[a,b]$ (for k=0,1,..,n-1) and their corresponding
causal Green's functions $G(x,y)=T(x,y)\theta(x-y)$ (given by
(\ref{maryam})), on subspace of ${\mc}^{\infty}[a,b]$ consisting of
functions which themselves and their derivatives to all orders are
zero at $x=a$, creates non-abelian group with operators
multiplication.}
\end{itemize}

Beside all differential operators
$\mathcal{O}(x,\rd_{x})=\big({\rd^{\,n}_{x}}+\sum^{n-1}_{k=0}
P_{k}(x) {\rd^{k}_{x}}\big)$ and their causal Green's functions
$G(x,y)=T(x,y)\theta(x-y)$, the group also contains
integro-differential operators and their inverses, coming from
mixing these two groups of operators. For example
$\mathcal{O}_{1}(x,\rd_{x}).T_{2}(x,y)\theta(x-y)$, acting on
$\phi(x)$ in the function space as
$\mathcal{O}_{1}(x,\rd_{x})(\int_{a}^{x}\, dz \,T_{2}(x,z)\phi(z))$,
and its inverse $\mathcal{O}_{2}(x,\rd_{x}).T_{1}(x,y)\theta(x-y)$
(where $\mathcal{O}^{-1}_{1}(x,\rd_{x})=T_{1}(x,y)\theta(x-y)$ and
$\mathcal{O}^{-1}_{2}(x,\rd_{x})=T_{2}(x,y)\theta(x-y)$).

\section{Green's function for linear differential
operators in one variable with other boundary conditions}
 \label{tecneq2}
In previous part the causal Green's function was derived according
to coefficient functions of linear differential operator, here in
this part it is shown that Green's function of general linear
differential operator, for other boundary conditions on $[a,b]$
(e.g. Sturm-Liouville problem), can also be derived according to
coefficient functions of differential operator. First we note that
homogeneous linear differential equation of degree $n$ in one
variable, $\big({\rd^{\,n}_{x}}+\sum^{n-1}_{k=0} P_{k}(x)
{\rd^{k}_{x}}\big)\,u(x)=0,$ for initial condition $\rd^{\,i}_{x}
\,u(a)=c_{i}$ for $ i=0,1,..,n-1$ can be converted to Volterra's
integral equation of second kind as; $\big(\mu(x)+\int_{a}^{x} dz\,
K(x,z) \,\mu(z)\big)=S(x),$ where $
K(x,z)=\big(\sum_{k=0}^{n-1}P_{k}(x){{(x-z)^{n-k-1}}\over(n-k-1)!}\big)$,
$u(x)=D(x)+\int_{a}^{x} dz {(x-z)^{n-1}\over (n-1)!} \mu(z)$, $\,
D(x)=\sum_{k=0}^{n-1}{c_{k}{(x-a)^{k}}\over k!}$ and $ S(x)=
-\sum_{i=0}^{n-1}\sum^{i}_{k=0}
c_{i}{{P_{k}(x)(x-a)^{i-k}}\over(i-k)!}$. Therefore we state the
following theorem;

\begin{theorem}{The Solution of homogeneous linear differential equation
 $\big({\rd^{\,n}_{x}}$ $+\sum^{n-1}_{k=0} P_{k}(x)\,{\rd^{k}_{x}}
 \big)\,u(x)=0,$ where $P_{i}(x)$ \,(\,$i=0,1,$ $...,n-1$) are in $\mc[a,b]$,
  with initial condition; $\rd^{\,i}_{x}$ $u(a)= c_{i}$ for $i=0,1,$ $..,n-1$, is given
  by;
 \be u(x)=D(x)+\int^{x}_{a}dz\,T(x,z) S(z),\label{el}\ee
where $\displaystyle\, D(x)=\sum_{k=0}^{n-1}{c_{k}{(x-a)^{k}}\over
k!}$,\,$\,\,\,\displaystyle S(x)= -\sum_{i=0}^{n-1}\sum^{i}_{k=0}
c_{i}{{P_{k}(x)(x-a)^{i-k}}\over(i-k)!}$ and $T(x,y)$ is given by
(\ref{hasan}).}
\end{theorem}

By using the above theorem for solutions $u_{r}(x)$ ($r=0,...,n-1$)
with initial conditions $\rd_{x}^{\,i}u_{r}(a)=\delta_{r,i}$ for
$i=0,...,n-1$, one can find $n$ linearly independent solutions;

\begin{itemize}{\item{\bf Fundamental solutions of homogenous linear ordinary differential equation.\vspace{0.2cm}}} {\em
$\displaystyle u_{r}(x)={{(x-a)^{r}}\over r!}+\int^{x}_{a}dz\,T(x,z)
(-\sum_{k=0}^{r} {{P_{k}(z)(z-a)^{r-k}}\over(r-k)!})$ for
$r=0,...,n-1$ are $n$ linearly independent solutions of
$\vspace{0.2cm}\displaystyle \big({\rd^{\,n}_{x}}$
$+\sum^{n-1}_{k=0} P_{k}(x)\,{\rd^{k}_{x}}
 \big)\,u(x)=0$. The Wornskian  is given by Abel's identity as;
$\displaystyle W(u_0,u_1,...u_{n-1})=e^{-\int^{x}_{a} dz
P_{n-1}(z)}$.}
\end{itemize}

It is easy to check that for boundary condition; $\rd^{\,i}_{x}u(b)=
c_{i}$ for $i=0,1,..,n-1$ the answer for homogenous differential
equation in $\mc[a,b]$ is given by (\ref{el}) in which $a$ is
replaced by $b$. One can derive solutions of homogenous differential
equation, with specific boundary conditions on either $a$ or $b$, by
using (\ref{el}). By substituting solutions of homogenous
differential equation with suitable boundary conditions (derived by
(\ref{el})) in to the expressions for Green's functions given by
method of variation of parameters, one can derive the Green's
function according to coefficient functions of differential
operator. For example the Green's function for Sturm-Lioville
problem; \be (\rd_{x}^{2}-P(x)) y(x)= g(x), \hspace{1cm}
\text{B.C.}\hspace{0.6cm}y(a)=y(b)=0, \ee by method of variation of
parameters is given by; $G(x,y)=(W(y))^{-1}$ $\big(u_{1}(x)u_2(y)$
$\theta(y-x)$ $+u_{2}(x)u_{1}(y)\,$ $\theta(x-y)\big)$, where
$u_{1}(x)$ and $u_{2}(x)$ are answers of corresponding homogenous
differential equation, $W(y)$ is the Wronskian and we have
$u_{1}(a)=u_{2}(b)=0$. If we take B.C. $u_{1}(a)=0$ and $\rd_x
u_{1}(a)=1$ then according to (\ref{el}); we have
$u_{1}(x)=(x-a)+\int^{x}_{a}dz T_{s}(x,z)P(z)(z-a)=T_{s}(x,a)$ (The
function $T_{s}(x,y)$ is given by (\ref{shahla})). On the other hand
by taking B.C. $u_{2}(b)=0$ and $\rd_x u_{2}(b)=1$ we have
$u_{2}(x)=T_{s}(x,b)$. The Wronskian of $u_{1}(x)$ and $u_{2}(x)$ is
constant and can be calculated easily at point $x=a$ as;
$W(y)=-T_{s}(a,b)$ therefore; \be G(x,y)=(-T_{s}(a,b))^{-1}{
\Big(T_{s}(x,a)T_{s}(y,b)\theta(y-x)+(x\leftrightarrow y) \Big)
}.\ee
\section{Conclusion }\label{disc}

Causal Green's function for general linear differential operator in
one variable was given by [theorem \ref{FinalGreen's}] as series of
integrals. Multiplicative property of causal Green's functions is
shown by [theorem \ref{adel}]. For differential operators which are
equal to multiplications of first order linear differential
operators, explicit formula (\ref{iraj1}), was given for causal
Green's functions. An infinite non-abelian group of operators on a
subspace of ${\mc}^{\infty}[a,b]$ is introduced. Analytic formula
for fundamental solutions of homogenous linear differential equation
in one variable was given via equation (\ref{kasai}) and for
specific boundary condition via equation (\ref{el}) as series of
integrals. By using equation (\ref{el}) and the method of variation
of parameters a way to derive Green's functions with arbitrary
boundary conditions in one variable according to coefficient
functions of differential operators was given.

\section*{Acknowledgements} I acknowledge {\em Farhang Loran} for
his helps and useful discussions. I dedicate this paper to memory of
my father {\em Iraj Kassaian}.


\end{document}